\documentclass[11pt]{article}

\usepackage[centertags]{amsmath}
\usepackage{amsfonts}
\usepackage{amssymb}
\usepackage{theorem}
\usepackage{epsfig}
\usepackage{color}
\usepackage{graphicx}
\usepackage{mathrsfs}
\usepackage[boldsans]{ccfonts}
\theoremstyle{break} \newtheorem{theorem}{Theorem}[section]
\theoremstyle{break} 
\theoremstyle{break} \newtheorem{definition}{Definition}[section]       
\theoremstyle{break} \newtheorem{lemma}[theorem]{Lemma}
\theoremstyle{break} 
\theoremstyle{break} 
\theoremstyle{break} 
\theoremstyle{break}
{\theorembodyfont{\rmfamily}\newtheorem{remark}[theorem]{Remark}}
{\theorembodyfont{\rmfamily}}
\theoremstyle{break} 
\theoremstyle{break} 
\theoremstyle{break} 
\theoremstyle{break} 
\numberwithin{equation}{section}

\newcommand{\D}{{\mathbb{D}}}
\newcommand{\C}{{\mathbb{C}}}

\def\Re{\mathop{{\rm Re}}}

\def\arcsinh{\mathop{{\rm arcsinh}}}

\newcommand{\hide}[1]{}

\begin{document}

\renewcommand{\thefootnote}{}
\stepcounter{footnote}
\begin{center}
{\bf \Large The theorems of Schottky and Landau for  analytic\\[2mm] functions
omitting the $\mathbf{n}$--th roots of unity
\footnote{2000 Mathematics Subject
 Classification: Primary  30E25, 35J65, 53A30}}
\end{center}
\renewcommand{\thefootnote}{}
\setcounter{footnote}{0}
\begin{center}
{\large Daniela Kraus and Oliver Roth}\\[4mm]

{\small To Stephan Ruscheweyh on the occassion of his 70th birthday}
\end{center}


\bigskip
\renewcommand{\thefootnote}{\arabic{footnote}}
\smallskip
\begin{center}
\begin{minipage}{13cm} {\bf Abstract.}   
We prove sharp  Landau-- and Schottky--type theorems for analytic
functions which omit the $n$--th roots of unity. The proofs are based on a
sharp lower bound for the Poincar\'e metric of the complex plane punctured at the
roots of unity.
\end{minipage}
\end{center}
\footnotetext{Research supported
by a DFG grant (RO 3462/3--2).}

\section{Introduction}

Let $\D:=\{z \in \C \, : \, |z|<1\}$ denote the open unit disk in the complex
plane $\C$.
In \cite{Hem79} J.~Hempel established sharp versions of the classical theorems
of Landau \cite{Lan04} and Schottky \cite{Sch04}, which are concerned with 
functions analytic  on $\D$ that omit the values $0$ and $1$. 
Related results have been obtained e.g.~by Ahlfors \cite{Ahl38}, Hayman
\cite{Hay} and Jenkins \cite{Jen}.

\medskip

In \cite{KRS} Hempel's work has been extended
 to meromorphic functions belonging to the classes ${\cal M}_{j,k,l}$, which are defined
as follows. Let $j,k,l \ge 2$ be integers (or $=\infty$)
such that
\begin{equation*} \label{eq:1}
 \frac{1}{j}+\frac{1}{k}+\frac{1}{l}<1 \, , 
\end{equation*}
and, for any domain $G \subseteq \C$, let
\begin{equation*} \label{eq:eq}
\hspace*{-0.4cm}\begin{array}{rl}
{\cal M}_{j,k,l}(G):=\{f \text{ meromorphic in } G \text{ s.t.} & \text{\,\,(i) all
  zeros of } f \text{ have
 order } \ge j,\\ &  \text{\,(ii) all} \text{ zeros}   \text{ of } f-1  \text{ have order } \ge k,\\
& \hspace*{0.6cm} \text{ and } \\ & \text{(iii) all
poles of } f \text{ have order }\ge l\}\, .
\end{array}
\end{equation*}
Note that ${\cal M}_{\infty,\infty,\infty}(\D)$ consists of all functions
analytic in $\D$ that omit $0$ and $1$. The motivation for introducing the classes
${\cal M}_{j,k,l}$ comes mainly from  Nevanlinna's celebrated Second
Fundamental Theorem which in particular implies
 that ${\cal M}_{j,k,l}(\C)$ contains only constant functions.
The corresponding normality result has been established by 
 Drasin \cite{Dra69}, who proved that ${\cal M}_{j,k,l}(\D)$ is a normal
 family,  thereby extending an earlier result of  Montel \cite[p.~125--126]{Mon27}.
Recently, sharp qualitative Landau-- and Schottky--type theorems for the
classes ${\cal M}_{j,k,l}$ have been obtained in \cite{KRS}.

 \medskip

In the present paper we employ the methods and results of \cite{KRS}  to
prove precise Landau-- and Schottky--type results for functions analytic in
$\D$ which omit the $n$--th roots of unity.
For any integer $n \ge 1$ we denote by
\begin{equation*} 
S_n=:\{e^{2\pi i\, j/n} \, : \, j=1, \ldots , n\}
\end{equation*} 
the set of $n$--th roots of unity. We say a  function $f$ analytic in $\D$
omits a set $S \subseteq \C$ if $f(\D) \cap S=\emptyset$. Further, it will be
convenient  to make the following formal definition.

\begin{definition}
For each integer $n\ge 2$ let 
\begin{equation*} \label{eq:Cn}
\gamma_n =\frac{1}{\lambda_{\C\backslash S_n}(e^{i\pi/n})} \, ,
\end{equation*}
where $\lambda_{\C \backslash S_n}$ denotes the density of the hyperbolic 
metric on $\C \backslash S_n$ with constant curvature $-1$.
\end{definition}
For details concerning the hyperbolic metric $\lambda_{\C \backslash S_n}(z)\,
|dz|$ we refer to Section \ref{sec:2}.

\medskip

Now we can state a sharp Landau--type theorem for functions which omit $S_n$ for
 a fixed integer $n \ge 2$. 

\begin{theorem} \label{thm:1}
Fix an integer $n \ge 2$.
Let $f(z)=a_0+a_1 z+\ldots $ be an analytic function on $\D$ and suppose that $f$  omits $S_n$. Then
\begin{equation} \label{eq:landau}
|a_1| \le \begin{cases}
2 |a_0| \sinh  \big( \arcsinh \left(\gamma_n\right)-\log|a_0|  \big) &
\quad\quad \text{ if } \, \,
|a_0| \le 1, \\[2mm]
2 |a_0| \big(  \gamma_n+ \log|a_0| \big) & \quad \quad \text{ if } \, \,|a_0| >1
\, . 
\end{cases}
\end{equation}
 Equality occurs in (\ref{eq:landau})   if and only if $f$ is a universal
covering map onto $\C \backslash S_n$ with $f(0)^n=\nolinebreak
-1$.
\end{theorem}

Here, for $a_0=0$ inequality (\ref{eq:landau}) is interpreted as
\begin{equation*}
|a_1|\le \lim_{a_0\to 0}\Big[ 2 |a_0| \sinh  \big( \arcsinh
\left(\gamma_n\right)-\log|a_0|  \big) \Big]= \frac{1 + 2 \gamma_n^2 + 
 2 \gamma_n \sqrt{1 + \gamma_n^2}}{\gamma_n + \sqrt{1 + \gamma_n^2}}\, .
\end{equation*}

\medskip

Some information about the constants $\gamma_n$ can be gleaned.
For this, recall that the standard hypergeometric function
$_2F_1(a,b,c;z)$ is defined as
$$ _2F_1(a,b,c;z)=\sum \limits_{k=0}^{\infty} \frac{(a)_k (b)_k}{(c)_k} z^k \,
, \qquad z \in \D\, , $$
where $(a)_k:=a (a+1) \cdots (a+k-1)$ is the Pochhammer symbol.

\begin{theorem}\label{thm:gamma_n}
\begin{itemize}
\item[(a)] For every integer $n\ge 2$, the constant $\gamma_n$ can be computed:
\begin{equation}\label{eq:value-gamma_n}
\hspace*{-1cm} \begin{array}{rcl} 
\gamma_n & =&   
\displaystyle \frac{ 2^{\frac{1}{n}}  \sin \left(\frac{\pi }{n}\right) \Gamma \left(\frac{1}{n}\right) \Gamma
   \left(\frac{n-1}{2 n}\right)^2}{n\, \pi} \textstyle
\, _2F_1\left(\frac{n-1}{2 n},\frac{n-1}{2 n},1;\frac{1}{2}\right) \, _2F_1\left(\frac{n-1}{2
   n},\frac{n-1}{2 n},\frac{n-1}{n};\frac{1}{2}\right) \\[4mm]
& & \displaystyle - \, \frac{  2^{\frac{1}{n}}  \pi \tan \left(\frac{\pi }{2
   n}\right)}{n} \,  \textstyle _2F_1\left(\frac{n-1}{2 n},\frac{n-1}{2 n},1;\frac{1}{2}\right)^2  \, .
\end{array}
\end{equation}

\item[(b)]
The constants $\gamma_n$ are strictly decreasing positive numbers, i.e.,
$\gamma_n>\gamma_{n+1}>0$ for $n=2,3,  \ldots$ and $\lim\limits_{n \to \infty}
\gamma_n=0$.

\item[(c)] For every integer $n\ge 2$: 
\begin{equation*}
\gamma_n =\dfrac{1}{\min \limits_{|z|=1} \lambda_{\C \backslash S_n}(z)}\,.
\end{equation*}
\end{itemize}
\end{theorem}

Here is a little overview of the first  values of the $\gamma_n$'s.

\begin{center}
\begin{tabular}{|@{\hspace{5mm}}c@{\hspace{5mm}}  |@{\hspace{5mm}}c@{\hspace{5mm}} |}
\hline 
 {$n$}  &  {$\gamma_n$}  \\\hline
 2 & 3.52993 \\\hline
 3& 1.79372\\\hline
4& 1.22801\\\hline 
5& 0.942245 \\\hline 
10& 0.445789 \\ \hline
100 & 0.0437768\\ \hline
1000& 0.00437689\\ \hline
\end{tabular}
\end{center}

\medskip

It is perhaps worth making some remarks. 

\begin{remark}\label{rem:0}
Since $(\gamma_n)$ is a strictly decreasing sequence, the
upper bound (\ref{eq:landau}) for $|a_1|$ provided by Theorem \ref{thm:1} is
strictly decreasing 
with respect to $n$. At first sight, this is a surprising fact, because
 a function which omits $S_{n+1}$ does not necessarily omit $S_{n}$.
\end{remark}

\begin{remark}[The case $\mathbf{n=2}$]\label{rem:landau}
The case $n=2$ of Theorem \ref{thm:1} is related to Hempel's sharp version of
Landau's theorem \cite{Hem79}, which says that
$$ |b_1| \le 2 |b_0| \left( \big|\log|b_0|\big|+\frac{\Gamma
    (\frac{1}{4})^4}{4 \pi^2} \right)$$
for every analytic function $g(z)=b_0+b_1 z+\ldots$ from $\D$ into $\C
\backslash \{0,1\}$. This is equivalent to
\begin{equation} \label{eq:hem1}
|a_1| \le 2 |a_0+1| \left( \left| \log \left|\frac{a_0+1}{2}\right| \right|+\frac{\Gamma
    (\frac{1}{4})^4}{4 \pi^2} \right)
\end{equation}
for every analytic function $f(z)=a_0+a_1 z+\ldots$ from $\D$ into $\C
\backslash S_2$. It can be shown that the estimate (\ref{eq:landau}) of
Theorem \ref{thm:1} is sharper than (\ref{eq:hem1}) for $a_0$ in a
neighborhood of the origin or in the right halfplane. On the other hand,
(\ref{eq:hem1}) is better than (\ref{eq:landau}) if $a_0$ is close to $-1$,
see Figure \ref{fig:1}.\\
\begin{figure}[h]
\centerline{\includegraphics[width=7cm]{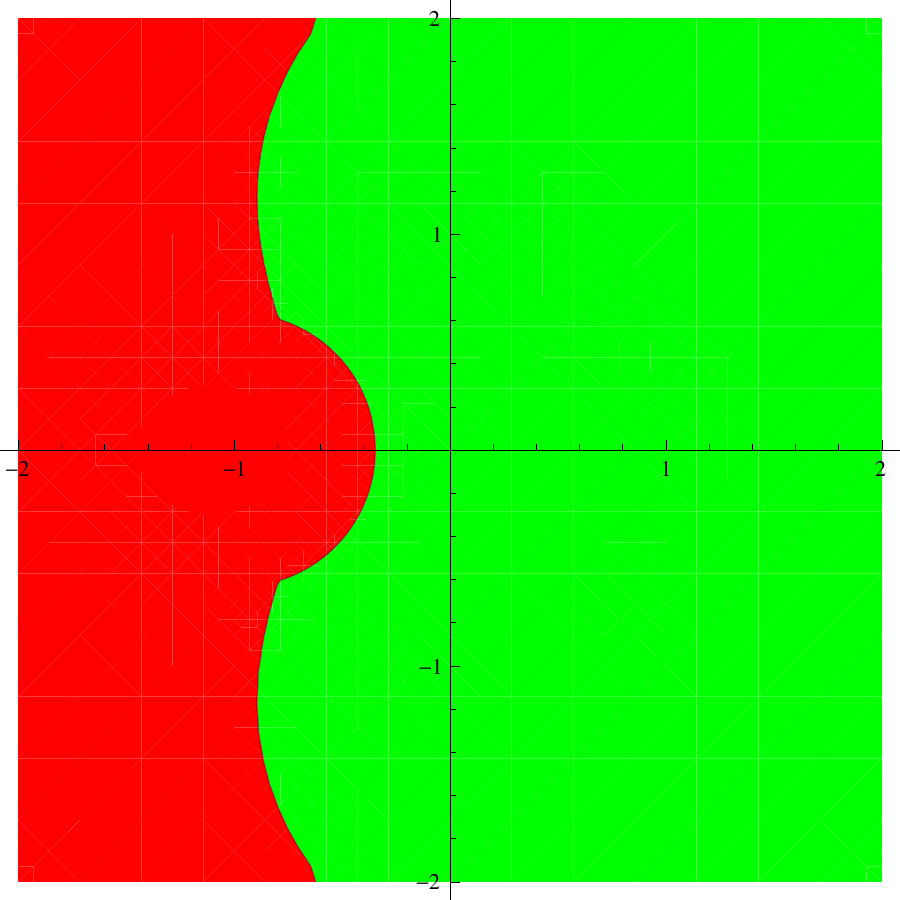} \hspace*{1cm} \includegraphics[width=7cm]{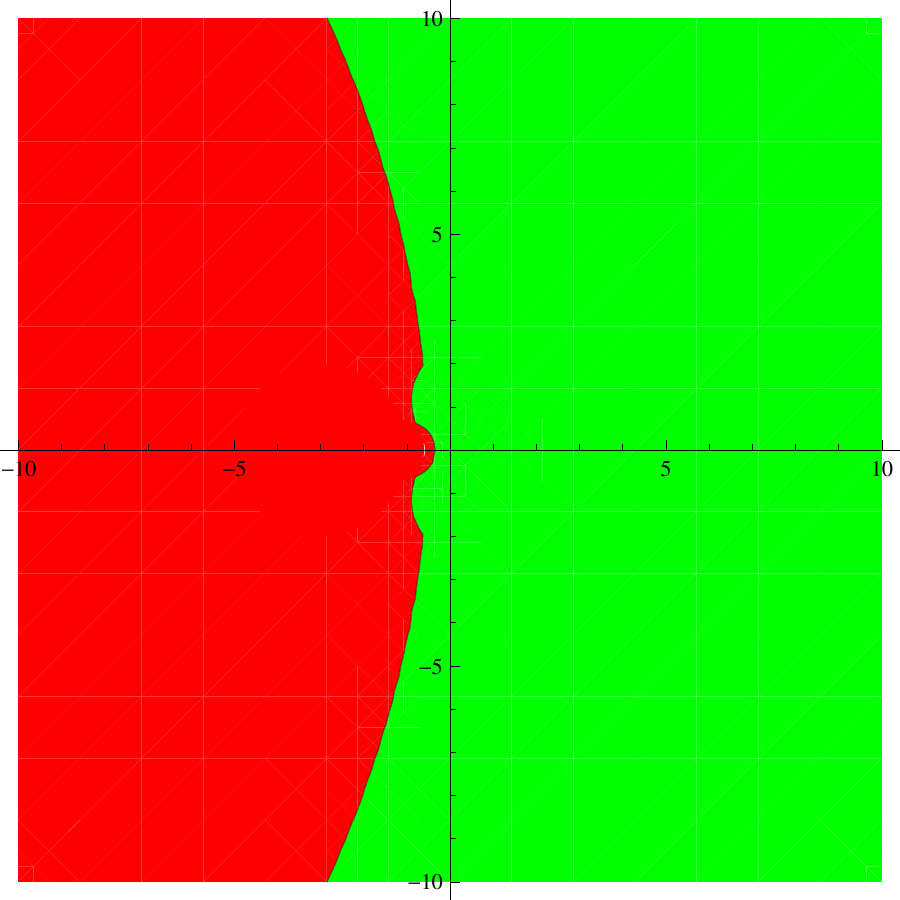}} 
\caption{The estimate (\ref{eq:landau}) is sharper than (\ref{eq:hem1}) for
  $a_0$ in the green region and the estimate
  (\ref{eq:hem1}) is sharper than (\ref{eq:landau}) in the red region.
} \label{fig:1}
\end{figure} 
\end{remark}

\begin{remark}[The case $\mathbf{n \to \infty}$]\label{rem:schwarz-pick}
In the limit $n \to \infty$, Theorem \ref{thm:1} reduces to the well--known
Schwarz--Pick lemma. In order to see this, note that in view of Theorem
\ref{thm:gamma_n} (b)
\begin{equation*}
 \lim \limits_{n \to \infty} \gamma_n=0 \, .
\end{equation*}
Hence, if a function $f(z)=a_0+a_1 z+\ldots$ analytic in $\D$ omits $\partial
\D=\{z \in \C: |z|=1\}$, then letting $n \to \infty$ in estimate (\ref{eq:landau}) yields
\begin{equation*} 
|a_1| \le \begin{cases}
\, 2 |a_0| \sinh  \left(   -\log|a_0|  \right)=1-|a_0|^2 & \quad\qquad
\text{ if }\, \, 
|a_0| \le 1, \\[2mm]
\, 2 |a_0| \log|a_0| & \quad \qquad \text{ if }\, \,  |a_0| >1
\, . 
\end{cases}
\end{equation*}
In particular, if $f$ maps $\D$ into $\D$, this implies the Schwarz--Pick
inequality $|a_1| \le 1-|a_0|^2$.
\end{remark}

We wish to emphasize that Theorem \ref{thm:1} might be viewed as
interpolation  between Landau's theorem ($n=2$; Remark \ref{rem:landau}) and the
Schwarz--Pick lemma ($n=\infty$; Remark \ref{rem:schwarz-pick}).

\medskip

The next result provides a sharp Schottky--type bound for analytic functions
in $\D$ which omit $S_n$ for fixed integer $n\ge 2$.
As usual $ \log^+x:=\max\{0,\log x\}$ for every positive real number $x$.

\begin{theorem} \label{thm:2}
Fix an integer $n \ge 2$.
Let $f$ be an analytic function on $\D$ and suppose that
$f$  omits $S_n$. Then for all $z \in \D$
\begin{equation} \label{eq:schottky}
\log |f(z)| \le \big( \gamma_n+\log^+|f(0)| \big)
\frac{1+|z|}{1-|z|}-\gamma_n\, .
\end{equation}
This estimate is sharp in the following sense: if $M>0$
is a constant such that
\begin{equation} \label{eq:schottky2}
\log |f(z)| \le \left[ M+\log^+|f(0)| \right]
\frac{1+|z|}{1-|z|}-M
\end{equation}
holds for all $z \in \D$ and all
analytic functions $f : \D \to \C \backslash S_n$, then $M \ge \gamma_n$.
\end{theorem}

Estimate (\ref{eq:schottky}) is rather crude at $z=0$ if $|f(0)| <1$, e.g.~if
$f(0)=0$. However, in this case one can combine the upper bound
(\ref{eq:schottky}) with Schwarz's lemma in order to arrive at the following result.

\begin{theorem}[Schwarz lemma for functions omitting the $\mathbf{n}$--th
  roots of unity] \label{thm:2a}
Let $n \ge 2$ be an integer and define
\begin{equation} \label{eq:Rn}
\begin{array}{rcl} 
R_n & = & \displaystyle 1+\gamma_n-\sqrt{\gamma_n^2+2 \gamma_n} \, .
\end{array}
\end{equation}

Suppose $f$ is an analytic function on $\D$ that  omits $S_n$. If $f(0)=0$, then
\begin{equation}  \label{eq:schwarz}
|f(z)| \le \frac{1}{R_n} \exp \left( \sqrt{\gamma_n^2+2 \gamma_n}-\gamma_n
\right)\, |z| \, , \qquad |z|<R_n \, . 
\end{equation}
\end{theorem}

Theorem \ref{thm:2a} merits some comment. 
\begin{remark} \label{rem:2}
We note (see Section \ref{sec:3.3}) that $(R_n)$ is a strictly increasing sequence with limit $1$ and
\begin{equation*}
\left( \frac{1}{R_n} \exp \left(  \sqrt{\gamma_n^2+2 \gamma_n}-\gamma_n
\right) \right) 
\end{equation*}
forms a strictly decreasing sequence, which converges to $1$ as $n$ tends to $\infty$.
Hence, letting $n \to \infty$   Theorem \ref{thm:2a}  reduces to Schwarz's
lemma. Consequently, the estimate (\ref{eq:schwarz}) is asymptotically sharp.  

\smallskip

A few examples for $R_n$ and $\frac{1}{R_n} \exp \left(  \sqrt{\gamma_n^2+2 \gamma_n}-\gamma_n
\right)$ are given in the next table.

\smallskip

\begin{center}
\begin{tabular}{|@{\hspace{5mm}}c@{\hspace{5mm}}
    |@{\hspace{5mm}}c@{\hspace{5mm}} | @{\hspace{5mm}}c@{\hspace{5mm}} |}
\hline 
 $ n$  &   $R_n$ & $\frac{1}{R_n} \exp \left(  \sqrt{\gamma_n^2+2 \gamma_n}-\gamma_n
\right)$  \\\hline
 2 & 0.111756 &21.7516 \\\hline
 3& 0.185105 &12.2035 \\\hline
4&   0.237023 & 9.0483  \\\hline 
5& 0.277218  & 7.43155  \\\hline 
10&0.401612  &  4.5297 \\\hline 
100 &0.744661 & 1.73354\\\hline
1000&0.910713 & 1.20059\\\hline
\end{tabular} 
\end{center}

We further remark that
\begin{equation*}
 |f'(0)| \le \frac{1}{R_n} \exp \left(  \sqrt{\gamma_n^2+2 \gamma_n}-\gamma_n\
\right)=1+8 \frac{\Gamma(5/4)}{\Gamma(3/4)} \frac{1}{\sqrt{n}}+ O(1/n) \qquad
(n \to \infty)  
\end{equation*}
for each $f \in \mathcal{A}_n:=\{f  : \D \to \C \backslash S_n
\text{ holomorphic} \, | \, f(0)=0\}$.
However,  if $f_n \in \mathcal{A}_n$  is a universal covering map 
of $\C \backslash S_n$, then we get 
\begin{equation*}
 \sup_{f \in \mathcal{A}_n} |f'(0)|=|f_n'(0)|
=\frac{\Gamma \left(\frac{1}{2} \left(1-\frac{1}{n}\right)\right)^2 \Gamma
   \left(\frac{1}{n}\right)}{n \Gamma \left(1-\frac{1}{n}\right) \Gamma \left(\frac{1}{2}
   \left(1+\frac{1}{n}\right)\right)^2}=
1+\frac{4 \log 2}{n}+O(1/n^2)  \qquad
(n \to \infty)  \,.
\end{equation*}
This shows that  estimate (\ref{eq:schwarz}) is asymptotically
sharp only ``up to order zero''.
\end{remark}

At the end of this introduction we wish to mention that a major motivation for
considering holomorphic functions which omit the $n$--th roots of unity has been the beautiful proof
of Montel's fundamental normality criterion based on the Zalcman lemma
(see \cite{Zal98}) where holomorphic functions which omit the $n$--th roots of
unit play a key role.

\medskip

In the following section we prove a sharp lower bound for the hyperbolic
metric $\lambda_{\C\backslash S_n}(z)\, |dz|$, see Theorem \ref{thm:3}. This result
is a key ingredient in the proof of Theorem \ref{thm:1}, but it is of some
interest in its own right. The proofs of the results in the Introduction are deferred to Section \ref{sec:3}.

\section{An explicit  sharp lower bound for the hyperbolic metric of 
the complex plane punctured at the $\mathbf{n}$--roots of unity}\label{sec:2}

Let $ U \subseteq \C$ be an open set. The quantity $\lambda(z)\, |dz|$ is
called a conformal metric on $U$, if $\lambda$ is a strictly positive function on $U$ and $\lambda$ is twice continuously
differentiable on $U$, i.e.~$\lambda \in C^2(U)$. If  in addition $\lambda$ satisfies
\begin{equation*}
\frac{\Delta \log\lambda(z)}{\lambda(z)^2}=1 \quad \text{ for all } z \in U,
\end{equation*}  
where $\Delta$ is the standard Laplace Operator, then we say $\lambda(z)\, |dz|$
has constant curvature $-1$ on $U$. It is a well-known fact that a hyperbolic domain
$\Omega \subset \C$, i.e.~$\Omega$ has at least two distinct boundary points, carries a
unique maximal conformal metric with constant curvature $-1$ (see, for
instance, Theorem 13.2 in \cite{Hei62}). This metric is called the hyperbolic
metric for $\Omega$ and is denoted by
\begin{equation*}
\lambda_{\Omega}(z)\,|dz|\,.
\end{equation*}
In particular,  
\begin{equation}\label{eq:hyp-1}
\mu (z)\le \lambda_{\Omega}(z)\, , \qquad z \in \Omega\,,
\end{equation}
for every conformal metric $\mu(z)\, |dz|$ on $\Omega$ with constant curvature $-1$.

\smallskip

We now turn to a sharp lower bound for the hyperbolic metric $\lambda_{\C
  \backslash S_n}(z)\, |dz|$.

\begin{theorem} \label{thm:3}
Let $n \ge 2$ be an integer. Then for all $z \in \C \backslash S_n$,
\begin{equation} \label{eq:m1}
 \lambda_{\C\backslash S_n}(z) \ge 
\begin{cases} 
\displaystyle 
\frac{1}{|z| \sinh \big(\arcsinh \left(\gamma_n\right)-\log|z| \big)}
  \qquad \quad & \text{ if } \, 
|z| \le 1 \, ,\\[4mm]
\displaystyle 
\frac{1}{|z| \left( \displaystyle \gamma_n+ \log|z| \right)}    \qquad \quad
&\text{ if } \, |z|>1\, .
\end{cases}
\end{equation}
 Equality holds in (\ref{eq:m1}) if and only if $z^n=-1$. 
\end{theorem}

Here, it is understood that
\begin{equation*}
\lambda_{\C\backslash S_n}(0)\ge \lim_{z\to 0}  \frac{1}{|z| \sinh \big(\arcsinh \left(\gamma_n\right)-\log|z| \big)}=
\frac{2 (\gamma_n + \sqrt{1 + \gamma_n^2})}{1 + 2 \gamma_n^2 +  2 \gamma_n \sqrt{1 + \gamma_n^2}}\,.
\end{equation*}

The proof of Theorem \ref{thm:3} is based on the theory developed in
\cite{KRS}, which deals with conformal metrics with cusps and corners. So let
us put together some of the needed details.

\medskip

Let $G \subseteq \C$ be a domain. 
 Suppose first, $p \in G$ and $\lambda(z)\, |dz|$ is a conformal metric on
$G \backslash \{ p \}$. Then  we say $\lambda(z)\, |dz|$ has a corner of order $\alpha<1$ at $p$, if
\begin{equation*}
\log \lambda(z)  = - \alpha \log |z-p|+O(1) \qquad \text{ as } z \to p\, ,
\end{equation*}
and a cusp at $p$, if 
\begin{equation*}
\log \lambda(z)  =
- \log |z-p|-\log \left( -\log |z-p| \right)+O(1)\qquad \text{ as } z \to p\, .
  \end{equation*}

\smallskip

If, however, $ \Delta_R:=\{z \in \C:|z|>R\} \subseteq G$ for some $R\ge 0$  and $\lambda(z)\, |dz|$ is a conformal metric on
$G$, then  
$\lambda(z)\, |dz|$ has a corner of order $\alpha<1$ at $\infty$, provided
\begin{equation*}
\log \lambda(z)  = - (2-\alpha )\log |z|+O(1) \qquad \text{ as } z \to \infty\, ,
\end{equation*}
and a cusp at $\infty$, if 
\begin{equation*}
\log \lambda(z)  = -\log|z|-\log \left(\log |z| \right)+O(1)\qquad \text{ as } z \to \infty\,.
\end{equation*}

We say a conformal metric $\lambda(z)\, |dz|$ has a singularity of order $\alpha \le 1$ at $p \in
\C \cup \{ \infty\}$, if either $\lambda(z)\, |dz|$ has a corner of
order $\alpha<1$ at $p$  or a cusp at $p$, if $\alpha=1$. Note, that the only
isolated singularities of a conformal metric with constant
curvature $-1$ are corners and cusps (see \cite{Nit57,KR2007}).

\smallskip

In particular, the hyperbolic metric $\lambda_{\Omega}(z)\, |dz|$ for a hyperbolic domain
$\Omega\subset \C$ has a cusp at every isolated boundary point of $\Omega$ and if
$\Delta_R \subseteq \Omega$ for some $R\ge 0$, then $\lambda_{\Omega}(z)\, |dz|$
has also a cusp at $\infty$, see \cite[\S 18]{Hei62} and
 \cite{Min97}.

\smallskip

We now record
two lemmas about cusps and corners of conformal metrics, which will be needed
later. The first lemma gives some information about
the remainder function $O(1)$ at a cusp and the second result is a removable
singularity theorem.

\begin{lemma}[cf.~Theorem 3.5 in \cite{KR2007}]\label{lem:cusp}
Let $G \subseteq \C$ be a domain. 
\begin{itemize}
\item[(a)] If $p \in G$ and $\lambda(z)|dz|$ is a conformal metric  with constant curvature $-1$ on $G\backslash
\{p\}$ which has a cusp at $p$, then
\begin{equation*}\label{eq:cusp}
\lim_{z\to p} \Big(\log \lambda(z) + \log |z-p|+\log \left( -\log |z-p|
  \right)\Big)=0\,.
\end{equation*}
\item[(b)]
If $\Delta_R \subset G$ and $\lambda(z)|dz|$ is a conformal metric  with
constant curvature $-1$ on $G$ which has a cusp at $\infty$, then
\begin{equation*}
\lim_{z\to \infty} \Big(\log \lambda(z) + \log |z|+\log \left( \log |z|
  \right)\Big)=0\, .
\end{equation*}
\end{itemize} 
\end{lemma}

\begin{lemma}[cf.~Theorem 1.1 in \cite{KR2007}]\label{lem:corner}
Let $G \subseteq \C$ be a domain, $p \in G$ and suppose $\lambda(z)|dz|$ is a
conformal metric  with constant curvature $-1$ on $G \backslash \{p \}$ which
has a corner of order $\alpha=0$ at $p$. Then the function $\log \lambda$ has a
$C^2$--extension to $G$. In particular, $\lambda(z)\, |dz|$ is a conformal
metric with constant curvature $-1$ on $G$.
\end{lemma}

It is a well--known fact that for real parameters $\alpha_1,\alpha_2, \alpha_3
\in  (0,1]$ which fulfill $\alpha_1+\alpha_2+\alpha_3>2$ there exists a
unique conformal metric
\begin{equation*}
\lambda_{\alpha_1,\alpha_2,\alpha_3}(z)\, |dz|\, , \qquad z \in \C\backslash \{0,1\}\,,
\end{equation*} 
with constant curvature $-1$ on $\C\backslash \{0,1\}$ and a singularity of
oder $\alpha_1$ at $z=0$, $\alpha_2$ at $z=1$ and $\alpha_3$ at $z=\infty$
such that
\begin{equation}\label{eq:gen-hyperbolic}
\mu(z)\le \lambda_{\alpha_1,\alpha_2,\alpha_3}(z)\, , \quad z \in \C
\backslash \{0,1 \}\,,
\end{equation}
for every conformal metric $\mu(z)\, |dz|$ on $\C\backslash \{0,1 \}$ with
constant curvature $-1$ which has a singularity of order $\beta_1 \le
\alpha_1$ at $z=0$, $\beta_2\le \alpha_2$ at $z=1$ and $\beta_3\le \alpha_3$
at $z=\infty$, see \cite[\S 18--\S 22]{Hei62}. 
This metric $\lambda_{\alpha_1,\alpha_2,\alpha_3}(z)\, |dz|$ is also
called the generalized hyperbolic metric on $\C \backslash \{0,1\}$ with singularities of order
$\alpha_1$, $\alpha_2$ and $\alpha_3$ at $z=0$, $z=1$ and $z=\infty$ respectively.

\medskip

In the proof of Theorem \ref{thm:3}, the conformal metric 
\begin{equation*}
\lambda_{1-1/n,1,1}(w) \, |dw|
\, , \qquad w \in \C\backslash \{0,1\},
\end{equation*}
where $n\ge 2$ is an integer, 
provides a key link to $\lambda_{\C\backslash S_n}(w)\, |dw|$.
In order to relate $\lambda_{\C\backslash S_n}(z)\, |dz|$ with $\lambda_{1-1/n,1,1}(z)\,
|dz|$ we use the concept of the
pullback of a conformal metric.

\smallskip

Suppose $\lambda(w)\, |dw|$ is a conformal metric on a domain 
$D\subseteq \C$ and $w=f(z)$ is a nonconstant analytic function from a  domain
$G \subseteq \C$ into $D$, then
\begin{equation*}
f^*\lambda(z)\,|dz|:=\lambda(f(z))\,|f'(z)|\, |dz|
\end{equation*}
is a conformal metric on $G\backslash \{ z\in G: f'(z)=0\}$
and is called the pullback of $\lambda(w)\, |dw|$ under the map $f$. If in
addition $\lambda(z)\, |dz|$ has constant curvature $-1$ on $D$, then it is
easy to see that
$f^*\lambda(z)\, |dz|$ has constant curvature $-1$ on $G\backslash \{ z\in G:
f'(z)=0\}$. 

\smallskip

Now, roughly speaking, the hyperbolic metric $\lambda_{\C \backslash S_n}(w) \, |dw|$ is the
pullback of the generalized hyperbolic metric $\lambda_{1-1/n,1,1}(z) \, |dz|$ under the map
$w=z^n$. This is the content of our next theorem.

\begin{theorem} \label{thm:3a}
Let $n \ge 2$ be an integer. Then
$$ \lambda_{\C \backslash S_n}(z)=n |z|^{n-1} \lambda_{1-1/n,1,1}(z^n) \, $$
for every $z \in \C \backslash (S_n\cup \{0\})$ and
\begin{equation*}
\lambda_{\C\backslash S_n}(0)=\lim_{z\to 0} n|z|^{n-1} \lambda_{1-1/n,1,1}(z^n)\,.
\end{equation*}
\end{theorem}

{\bf Proof.} Let $f:\C \to \C$, $f(z)=z^n$. Note that $f$ maps $\C\backslash
(S_n \cup \{0 \})$ onto $\C \backslash \{0,1 \}$. Thus
\begin{equation*}
\lambda_n(z) \, |dz|:=f^*\lambda_{1-1/n,1,1}(z)\, |dz|=n |z|^{n-1}
\lambda_{1-1/n,1,1}(z^n) \, |dz| \, , \qquad z \in \C\backslash (S_n \cup \{0 \})\,,
\end{equation*}
defines a conformal metric on  $\C\backslash (S_n \cup \{0 \})$  with constant
curvature $-1$.

\smallskip

Our objective is to show that $\lambda_n(z)\, |dz|$ extends to be a conformal
metric with constant curvature $-1$ on all of $\C \backslash S_n$ and
\begin{equation*}
\lambda_n(z)\,
|dz|=\lambda_{\C \backslash S_n}(z)\, |dz|\quad \text{ on } \, \C\backslash S_n\,.
\end{equation*}

\smallskip

We first determine the type of singularities of $\lambda_n(z)\, |dz|$ at
the points $p \in S_n \cup \{0, \infty\}$.
As $\lambda_{1-1/n,1,1}(z) \, |dz|$ has a
cusp at $p=1$, it follows that $\lambda_n(z) \, |dz|$ has a cusp at every
point $p \in S_n$. Similarly, the cusp of $\lambda_{1-1/n,1,1}(z) \, |dz|$ at
$p=\infty$ results in a cusp of $\lambda_n(z)\, |dz|$ at $p=0$.
Now, since $\lambda_{1-1/n,1,1}(z) \, |dz|$ has a corner of order
$1-1/n$ at $p=0$, we see that $\lambda_n(z)\, |dz|$ is bounded at $p=0$,
i.e.~$\lambda_n(z)\, |dz|$ has a corner of order $\alpha=0$ at $p=0$. By
appeal to Lemma \ref{lem:corner}, $\lambda_n(z)\, |dz|$ is a conformal metric
on $\C\backslash S_n$ with constant curvature $-1$.

\smallskip

In summary, $\lambda_n(z)\, |dz|$ is a conformal metric
 with constant curvature $-1$ on $\C\backslash S_n$  which has
 a cusp at each point $p \in S_n \cup \{ \infty\}$. 

\smallskip

It follows that $\lambda_n(z)\le \lambda_{\C \backslash S_n}(z)$ for $z \in \C
\backslash S_n$, because $\lambda_{\C \backslash S_n}(z)\, |dz|$ is the
hyperbolic metric on $\C \backslash S_n$, see (\ref{eq:hyp-1}).
We next consider on $\C \backslash S_n$ the nonnegative function
\begin{equation*}
s(z)=\log \left(\frac{\lambda_{\C \backslash
   S_n}(z)}{\lambda_n(z)}\right)\, , \qquad z \in \C\backslash S_n\,.
\end{equation*}
Then, $s$ is subharmonic on $\C \backslash S_n$, since 
\begin{equation*}
\Delta s(z)=\Delta \log \lambda_{\C \backslash S_n}(z)-\Delta \log
\lambda_n(z)= \lambda_{\C\backslash S_n}(z)^2-\lambda_n(z)^2 \ge 0\,, \qquad z
\in \C \backslash S_n\, ,
\end{equation*}
which is an immediate consequence of the fact that $\lambda_n(z) \, |dz|$ and $\lambda_{\C \backslash S_n}(z) \,
|dz|$ have  constant curvature $-1$ on $\C \backslash S_n$.
Bearing in mind that both $\lambda_n(z)\, |dz|$ and $\lambda_{\C \backslash
 S_n}(z)\, |dz|$ have a cusp at each point $p\in S_n\cup \{\infty\}$, we conclude with the
  help of Lemma \ref{lem:cusp} that
\begin{equation*}
\lim_{z \to p} s(z)=0 \quad \text{ for each } \, p \in S_n \qquad \text{and}
  \qquad 
\lim_{z \to \infty} s(z)=0\,.
\end{equation*}
Thus, the maximum principle for subharmonic functions (see e.g.~Corollary 2.3.3 in \cite{Ran95}) implies $s(z) \le
0$ for all $z \in \C \backslash S_n$. Consequently, $s \equiv 0$ on $\C \backslash
S_n$. This gives the desired result
\begin{equation*}
\lambda_n(z) \, |dz|=\lambda_{\C \backslash S_n}(z) \, |dz| \, \text{ for } \, z \in \C
\backslash S_n \quad \text{ and } \quad \lambda_{\C \backslash S_n}(0)=\lim_{z\to 0} \lambda_n(z)\,.
\end{equation*}\\[-13mm]\phantom{m}
\hfill{$\blacksquare$}

\medskip

Now, Theorem \ref{thm:3a} enables us to prove Theorem \ref{thm:3} by applying
the following sharp lower bound for the generalized hyperbolic metric
$\lambda_{1-1/n,1,1}(z)\, |dz|$.

\begin{theorem}\label{thm:3b}
Let $n \ge 2$ be an integer. 
 Then for all $w \in \C \backslash \{0,1\}$,
\begin{equation*}
 \lambda_{1-1/n,1,1}(w) \ge 
\begin{cases}
\displaystyle \frac{1}{n |w| \sinh \left[ \arcsinh \left( \displaystyle
      \frac{1}{n \, \lambda_{1-1/n,1,1}(-1)} \right) - \frac{1}{n} \log|w|  \right]} &
\qquad \quad \text{ if }\,\, |w| \le 1 \, , \\
&  \phantom{  \text{ if } } \\
\displaystyle \frac{1}{|w| \left( \displaystyle \frac{1}{\lambda_{1-1/n,1,1}(-1)}+ \log|w|
  \right)}  & \qquad \quad \text{ if } \, \,  |w|>1\,,  
\end{cases}
\end{equation*}
 with equality if and only if
$w=-1$.
\end{theorem}

{\bf Proof.}
Theorem \ref{thm:3b} can be obtained from Theorem 2.2 in
\cite{KRS}. In fact, Theorem 2.2 in \cite{KRS} provides a sharp lower bound for the generalized
hyperbolic metric $\lambda_{\alpha_1,\alpha_2,\alpha_3}(z)\, |dz|$, when $\alpha_1,\alpha_2,\alpha_3 \in (0,1]$ and
$\alpha_1+\alpha_2+\alpha_3>2$. Choosing $\alpha_1=1-1/n$ and letting $\alpha_2 \nearrow
1$ and $\alpha_3 \nearrow 1$ gives Theorem \ref{thm:3b}. \hfill{$\blacksquare$}

\medskip

Before passing to the proof of Theorem \ref{thm:3} let us recall an important
result about conformal metrics.

\begin{lemma}[Theorem 7.1 in \cite{Hei62}]\label{lem:equality}
Let $G \subseteq \C$ be a domain, $\lambda(z)\, |dz|$ a conformal metric with
constant curvature $-1$ on $G$ and let $\mu$ be a nonnegative on $G$
function such that $\mu(z)\, |dz|$ is a conformal metric on $\{z\in G:
\mu(z)>0\}$ with constant curvature $-1$. Suppose $\mu(z)\le \lambda(z)$ for all $z \in G$. Then either
$\mu<\lambda$ on $G$ or else $\mu \equiv \lambda$ on $G$.
\end{lemma}

{\bf Proof of Theorem \ref{thm:3}.}
By Theorem \ref{thm:3a} we see that
\begin{equation}\label{eq:gamma_n-2}
\frac{1}{
n \lambda_{1-1/n,1,1}(-1)}=\frac{1}{\lambda_{\C\backslash
  S_n}(e^{i\pi/n})}=\gamma_n\,.
\end{equation}

Now combining Theorem \ref{thm:3a} and Theorem \ref{thm:3b} gives inequality
(\ref{eq:m1}).

\smallskip

It remains to show that equality occurs in (\ref{eq:m1}) if and only if
$z^n=-1$. This is clear if $z\not=0$ in view of Theorem \ref{thm:3a} and the
case of equality of Theorem \ref{thm:3b}. Therefore we need only show that for
$z=0$ strict inequality holds in (\ref{eq:m1}), i.e.,
\begin{equation*}
\lambda_{\C\backslash S_n}(0)>\lim_{z\to 0} 
\frac{1}{ |z| \sinh \left( \arcsinh \left( \displaystyle
      \gamma_n \right) -  \log|z|  \right)}
=\frac{2 (\gamma_n + \sqrt{1 + \gamma_n^2})}{1 + 2 \gamma_n^2 + 
 2 \gamma_n \sqrt{1 + \gamma_n^2}}\,.
\end{equation*}
For this we define $\mu$ on $\D$  by
\begin{equation*}
\mu(z)=\begin{cases}
\dfrac{2 (\gamma_n + \sqrt{1 + \gamma_n^2})}{1 + 2 \gamma_n^2 + 
 2 \gamma_n \sqrt{1 + \gamma_n^2}} \qquad \quad & \text{ if } \, \, z=0,\\[6mm]

\dfrac{1}{|z|\sinh\left(\arcsinh(\gamma_n)-\log|z| \right) } \qquad \quad &
  \text{ if } \, \, z\in \D\backslash \{0 \}\,.
\end{cases}
\end{equation*}
Since we have already observed that strict inequality occurs in (\ref{eq:m1}) 
for every $z \in \C$ with $z\not=0$ and $z^n=-1$, respectively, we have in particular
\begin{equation}\label{eq:mu}
\mu(z)< \lambda_{\C\backslash S_n}(z) \quad \text{ for all } z \in \D
\backslash \{0 \}\,.
\end{equation}
Now, it is easily checked that $\mu(z)\, |dz|$ is a conformal metric with constant
curvature $-1$ on $\D \backslash \{0 \}$. Since  $\mu$ is continuous at $z=0$
and $\mu(0)>0$, Lemma \ref{lem:corner} tells us that $\mu(z)\, |dz|$ is a
conformal metric with constant curvature $-1$ on the entire unit disk $\D$, so 
$\mu(z)\le \lambda_{\C\backslash S_n}(z)$ for all $z \in \D$ because of (\ref{eq:mu}). 
We are now in a position to apply Lemma \ref{lem:equality} for $G=\D$, and
obtain, using again (\ref{eq:mu}), that $\mu(z)<\lambda_{\C \backslash S_n}(z)$
for all $z \in \D$. This completes the proof.
\hfill{$\blacksquare$}

\section{Proofs}
\label{sec:3}

\subsection{Proofs of Theorem \ref{thm:1} and Theorem \ref{thm:gamma_n}}
Before proving Theorem \ref{thm:1} we digress to recall the Ahlfors--Schwarz
lemma.

\begin{lemma}[cf.~\cite{Ahl38} and \S 25 in \cite{Hei62}] \label{lem:ahlfors}
Let $G\subset \C$ be a hyperbolic domain,  suppose that $\lambda(z)\, |dz|$
is a conformal metric with constant curvature $-1$ on $G$, and let $f:\D \to G$
be a nonconstant analytic function. Then
\begin{equation}\label{eq:as0}
f^* \lambda(z)\le \lambda_{\D}(z)\, , \qquad z \in \D\,,
\end{equation}
where $\lambda_{\D}(z)\, |dz|=\dfrac{2}{1-|z|^2}\, |dz|$ is the hyperbolic
metric for $\D$.

\smallskip

Equality occurs in (\ref{eq:as0}) if and only if $f:\D \to G$ is
a universal covering map and $\lambda(z)\, |dz|=\lambda_{G}(z)\, |dz|$.
\end{lemma}
 
\medskip

{\bf Proof of  Theorem \ref{thm:1}.}  
Let $f : \D \to \C \backslash S_n$  be analytic and nonconstant. Then by Lemma \ref{lem:ahlfors}
\begin{equation} \label{eq:as1}
 f^*\lambda_{\C \backslash S_n}(z)=\lambda_{\C \backslash S_n}(f(z)) \,
|f'(z)| \le \frac{2}{1-|z|^2}\, , \qquad z \in \D \, .
\end{equation}

Thus (\ref{eq:as1}) together with Theorem \ref{thm:3} implies
\begin{equation} \label{eq:as2}
\frac{2}{1-|z|^2} \ge
\begin{cases}
\displaystyle \frac{|f'(z)|}{ |f(z)| \sinh \big(\arcsinh(\gamma_n) -\log|f(z)|  \big)} &
\quad \quad \text{ if } \, \, |f(z)| \le 1 \, ,  \\[6mm]
\displaystyle \frac{|f'(z)|}{|f(z)| \left( \gamma_n+ \log|f(z)|
  \right)}  & \quad \quad \text{ if } \, \,  |f(z)|>1  \, . 
\end{cases}
\end{equation}
As before  
\begin{equation*}
\frac{|f'(z)|}{ |f(z)| \sinh \big(\arcsinh(\gamma_n) -\log|f(z)|
  \big)}
\end{equation*}
 should be interpreted to be
\begin{equation*}
\lim_{f(z)\to 0} \frac{|f'(z)|}{ |f(z)| \sinh \big(\arcsinh(\gamma_n) -\log|f(z)|
  \big)}
\end{equation*}
 if $f(z)=0$.

\medskip

To obtain (\ref{eq:landau}) we just insert $z=0$ into (\ref{eq:as2}). 

\smallskip

We now turn to the case of equality. Note that equality holds in
(\ref{eq:landau}) if and only if equality occurs in (\ref{eq:as2}) for
$z=0$. In view of Theorem \ref{thm:3} and (\ref{eq:as1}) this is the case if
and only if 

\begin{equation*}
\lambda_{\C \backslash S_n}(f(0))|f'(0)|=\lambda_{\D}(0) \quad \text{ and } \quad f(0)^n=-1\,.
\end{equation*}

Hence, by Lemma \ref{lem:ahlfors} (case of equality) we see 
that equality occurs in (\ref{eq:landau}) if and only if  $f: \nolinebreak \D \to \C\backslash S_n$
is a universal covering with $f(0)^n=-1$.
\hfill{$\blacksquare$}

\bigskip

We next turn to the proof of Theorem \ref{thm:gamma_n}. 
In order to determine the value of the constant $\gamma_n$, we will make use of the following explicit
formula for $\lambda_{1-1/n,1,1}(z)\, |dz|$.  For this it is convenient to
define
\begin{equation*}
\varphi_1(z):=_2\!\!F_1\left (\frac{n-1}{2n},\frac{n-1}{2n},\frac{n-1}{n};z \right) \quad
\text{and} \quad \varphi_2(z):=_2\!\!F_1\left (\frac{n-1}{2n},
    \frac{n-1}{2n},1;1-z\right)\,.
\end{equation*}
Note that $\varphi_1$ has an analytic continuation to $\C\backslash [1,+\infty)$ and $\varphi_2$
has an analytic continuation to $\C\backslash (-\infty,0]$.

\begin{theorem}\label{thm:formula}
Let $n\ge 2$ be an integer. Then for $z \in \C \backslash\{0,1\}$

\begin{equation}\label{eq:formula}
\lambda_{1-1/n,1,1}(z)=\frac{1}{|z|^{1-1/n}|1-z| }\, \, 
\frac{1}{\frac{K_2}{2K_3}\, \,
    |\varphi_2(z)|^2+\frac{1}{K_3}\Re\Big(\varphi_1(z)
    \varphi_2(\overline{z})\Big)}\, , 
\end{equation}

where
\begin{equation*}
K_2:=-\frac{\Gamma\left(\frac{n+1}{2n}\right)^2}{\Gamma \left(\frac{1}{n}\right)} \quad \text{ and }
\quad K_3:= \frac{\Gamma \left( \frac{n-1}{n} \right)}{\Gamma \left(\frac{n-1}{2n}\right)^2}\,.
\end{equation*}
\end{theorem}

Here, for $x \in  (1,
\infty)$ the right--hand side of equation (\ref{eq:formula})  must be interpreted  to mean 
\begin{equation*}
\lim_{z\to x} \frac{1}{|z|^{1-1/n}|1-z| }\, \, 
\frac{1}{\frac{K_2}{2K_3}\, \,
    |\varphi_2(z)|^2+\frac{1}{K_3}\Re\Big(\varphi_1(z)
    \varphi_2(\overline{z})\Big)}\,. 
\end{equation*}

{\bf Proof.}
We note that Theorem 2.1 in \cite{KRS} displays a formula for the generalized
hyperbolic metric $\lambda_{\alpha_1,\alpha_2, \alpha_3}(z)\, |dz|$ for
$0<\alpha_1,\alpha_2<1$, $0<\alpha_3\le 1$ and $\alpha_1+\alpha_2+\alpha_3>2$ in terms of the functions $\varphi_1$ and
$\varphi_2$. Now choosing in this particular theorem  $\alpha_1=1-1/n$,
$\alpha_3=1$ and letting $\alpha_2\nearrow 1$ yields  (\ref{eq:formula}).
\hfill{$\blacksquare$}

\medskip

It is convenient to record different representations for
the constants $K_3$ and $K_2/K_3$ for later reference.

\smallskip

Using the well--known reflection formula 
\begin{equation}\label{eq:reflection}
\Gamma(z)\, \Gamma(1-z)=\frac{\pi}{\sin(\pi z)}\, ,\qquad z \in \C\backslash
\mathbb{Z}\, ,
\end{equation}
(cf.~\cite[p.~77]{AS84}), we get
\begin{equation}\label{eq:K3}
K_3=\frac{\pi}{\sin\left(\frac{\pi}{n}\right)}\,
  \frac{1}{\Gamma\left(\frac{1}{n}\right) \, \Gamma\left(\frac{n-1}{2n}\right)^2}
\end{equation}
and
\begin{equation}\label{eq:K2/K3}
\begin{split}
\frac{K_2}{K_3}&=- \frac{\Gamma \left(\frac{n+1}{2n} \right)^2}{\Gamma \left(\frac{1}{n}\right)}\,
\frac{\sin \left ( \frac{\pi}{n} \right)}{\pi} \, \, \Gamma \left(\frac{1}{n} \right) \, \, \Gamma
\left (\frac{n-1}{2n}\right)^2= -\frac{\sin ( \frac{\pi}{n})}{\pi} \,\, \Gamma
\left(\frac{n-1}{2n}\right)^2 \, \,\Gamma\left(\frac{n+1}{2n}\right)^2\\[2mm]
&=-\frac{\sin \left( \frac{\pi}{n}\right)}{\pi}\, \, \frac{\pi^2}{\left(\sin
    \left (
    \frac{\pi(n-1)}{2n}\right) \right)^2}= -2\pi \tan\left(\frac{\pi}{2n}\right)\,.
\end{split}
\end{equation}

{\bf Proof of Theorem \ref{thm:gamma_n}.} 

Recall that
\begin{equation*}
\gamma_n=\frac{1}{n \lambda_{1-1/n,1,1}(-1)}\,,
\end{equation*}
see (\ref{eq:gamma_n-2}).
\begin{itemize}

\item[(a)]
To compute $\gamma_n$ we consider the M\"obius transformation $T(z)=z/(z-1)$, which fixes $z=0$ and
interchanges $z=1$ with $z=\infty$.  Note that $\mu(z)\,
|dz|:=T^*\lambda_{1-1/n,1,1}(z)\, |dz|$ is a conformal metric on $\C
\backslash \{0,1\}$ with constant curvature $-1$  and a corner of order
$1-1/n$ at $z=0$ and a cusp at $z=1$ and $z = \infty$. Thus 
\begin{equation*}\label{eq:in2}
\mu(z)\le \lambda_{1-1/n,1,1}(z)\, , \qquad z \in \C \backslash \{0,1\}\,,
\end{equation*} 
by (\ref{eq:gen-hyperbolic}). Since $T\circ T= \text{Id}$, we observe that
$T^*\mu(z)\, |dz|=\lambda_{1-1/n,1,1}(z)\, |dz|$ and 
the latter inequality implies
\begin{equation*}
\lambda_{1-1/n,1,1}(z)\le \mu(z)\, , \qquad z\in \C \backslash \{0,1\}\,.
\end{equation*}
Consequently,
\begin{equation*}
\lambda_{1-1/n,1,1}(z)= \lambda_{1-1/n,1,1}(T(z)) |T'(z)|\, ,\qquad z \in
\C\backslash \{0,1\}\,.
\end{equation*}
Thus we get
\begin{equation*}
\lambda_{1-1/n,1,1}(-1)=\frac{\lambda_{1-1/n,1,1}(1/2)}{4} \, .
\end{equation*}
Inserting now $z=1/2$ into (\ref{eq:formula}) and using (\ref{eq:K3}) and
(\ref{eq:K2/K3}) gives
the claimed expression (\ref{eq:value-gamma_n}) for $\gamma_n$.

\item[(b)] 
By (\ref{eq:gen-hyperbolic}) we have, for all integers $m>n\ge 2$,
\begin{equation*}\label{eq:mono}
\lambda_{1-1/n,1,1}(w)\le\lambda_{1-1/m,1,1}(w)
\, , \qquad w \in \C \backslash \{0,1\}\,. 
\end{equation*}
This  shows that 
\begin{equation*}
0<\gamma_n=\frac{1}{n \lambda_{1-1/n,1,1}(-1)} < \frac{1}{m
  \lambda_{1-1/n,1,1}(-1)}=\gamma_m
\end{equation*}
for $m>n\ge 2$.

\smallskip

To verify  $\lim_{n \to \infty}\gamma_n=0$ we consider the sequence
\begin{equation*}
(n\gamma_n)=\left( \frac{1}{\lambda_{1-1/n,1,1}(-1)}\right)\,. 
\end{equation*}
Since 
\begin{equation*}
\lim_{n \to \infty} \lambda_{1-1/n,1,1}(-1)=\lambda_{1,1,1}(-1)=\lambda_{\C
  \backslash \{0,1 \}}(-1)
\end{equation*}
and
\begin{equation*}
\lambda_{\C
  \backslash \{0,1 \}}(-1)=\frac{4\pi^2}{\Gamma(1/4)^4}\, ,
\end{equation*}
see \cite{Hem79}, we get 
\begin{equation*}
\lim_{n \to \infty} n \gamma_n=\frac{\Gamma(1/4)^4}{4\pi^2}\quad \text{ and } \quad
\lim_{n \to \infty} \gamma_n=0\,.
\end{equation*}

\item[(c)]
It is a known fact, that
\begin{equation*}
\lambda_{1-1/n,1,1}(-1)=\min_{|z|=1}\lambda_{1-1/n,1,1}(z)\, ,
\end{equation*}
see, for instance, Theorem 4.10 in \cite{KRS}. Thus Theorem \ref{thm:3a} implies

\end{itemize}

\begin{equation*}
\min \limits_{|z|=1} \lambda_{\C\backslash S_n}(z)=
\min_{|z|=1}n\lambda_{1-1/n,1,1}(z^n)= n \lambda_{1-1/n,1,1}(-1)=\lambda_{\C\backslash S_n}(e^{i\pi/n})\,.
\end{equation*}\\[-10mm]\phantom{m}
\hfill{$\blacksquare$}

\subsection{Proof of Theorem \ref{thm:2}}

{\bf Proof of Theorem \ref{thm:2}.}
To prove (\ref{eq:schottky})
we proceed exactly as in the proof of Theorem 1 in \cite{LQ2007}
(see also the proof Theorem 2.8 in \cite{KRS}). Here are the details. Note
that (\ref{eq:schottky}) is true for $z=0$ and $|f(z)|<1$.
So pick a point $z \in \D\backslash \{0 \}$ such that $|f(z)|>1$.
Consider the curve $\gamma(t):=f(t \eta)$ for $t \in [0,|z|]$, where  
$\eta:=z/|z|$. If $\gamma$ lies completely outside $\overline{\D}$, then by (\ref{eq:as2})
\begin{equation}\label{eq:path}
\frac{|f'(t\eta)|}{|f(t\eta)| \left( \gamma_n+ \log|f(t\eta)|
  \right)} \le \frac{2}{1-t^2} \, , \qquad t \in [0,|z|] \, .
\end{equation}
Using  
\begin{equation*}
\frac{d}{dt}\left( |f(t\eta)| \right) \le |f'(t \eta)|\,,
\end{equation*}
 we obtain
\begin{equation*}
\log \left( \frac{\gamma_n +\log|f(z)|}{ \gamma_n +\log|f(0)|}\right)=
\int \limits_{0}^{|z|}\frac{\frac{d}{dt} \left(|f(t\eta)| \right)}{|f(t\eta)(\gamma_n +
  \log|f(t\eta)|)}\, dt \le
\int \limits_0^{|z|} \frac{2}{1-t^2}\, dt =\log\frac{1+|z|}{1-|z|} \,,
\end{equation*}
which leads to the desired  estimate (\ref{eq:schottky}). 

\smallskip 

If $\gamma$ does not completely lie outside $\overline{\D}$ we use an argument
similar to the one above. Let $\gamma(t^*)$ be the ``last'' point of
$\gamma(t)$ which lies inside $\overline{\D}$. Now integrating (\ref{eq:path})
just from $t^*$ to $|z|$ yields
\begin{equation*}
\log |f(z)| \le  \gamma_n \frac{1+|z|}{1-|z|}-\gamma_n  \, , 
\end{equation*}
which is stronger than (\ref{eq:schottky}). 

\medskip

Now suppose that (\ref{eq:schottky2}) holds for each analytic map $f$ from
$\D$ to $\C \backslash S_n$. 
Let $\varphi:\D \to \C\backslash S_n$ denote a universal covering map with
$\varphi(0)^n=-1$.  Then it follows from (\ref{eq:schottky2}) 
that
\begin{equation*}
\log|\varphi(z)|\le 2M \frac{|z|}{1-|z|} \qquad \text{ for all } z \in \D\,.
\end{equation*}

We further note that $|\varphi(0)|=1$ and hence there exists
a disk $\D_{r}:=\{z \in \C: |z|<r\}$, $0<r<1$, such that $\varphi(z)\not=0$ for all $z \in
\D_r$. 

\smallskip

Choose $\eta \in \partial \D$ such that
\begin{equation*}
\eta \frac{\varphi'(0)}{\varphi(0)}=\left| \frac{\varphi'(0)}{\varphi(0)} \right| =|\varphi'(0)|
\end{equation*}
and define 
\begin{equation*}
h(t)=\log|\varphi(\eta t)|\,, \qquad t\in (-r,r)\,.
\end{equation*}
Since $h(0)=0$ we conclude
\begin{equation*}
\frac{h(t)-h(0)}{t} \le 2M \frac{1}{1-t}\,.
\end{equation*}
Letting $t \to 0$, we obtain $h'(0)\le 2M$. On the other hand a simple
computation gives 
\begin{equation*}
h'(0)=\Re\left(\eta \frac{\varphi'(0)}{\varphi(0)}  \right)=|\varphi'(0)|\,.
\end{equation*}
 As $\varphi$ is a universal covering map of
$\C \backslash S_n$ with $\varphi(0)^n=-1$ Theorem \ref{thm:1} shows
$|\varphi'(0)|=2\gamma_n$. Putting these facts together, we see that $\gamma_n
\le M$.
\hfill{$\blacksquare$}

\subsection{ Proofs of Theorem \ref{thm:2a} and Remark \ref{rem:2}}
\label{sec:3.3}

{\bf Proof of Theorem \ref{thm:2a}.}
Choose $R \in (0,1)$.
Since $f(0)=0$,  estimate (\ref{eq:schottky}) implies 
\begin{equation*}
 |f(z)| \le \exp \left( 2 \gamma_n \frac{|z|}{1-|z|} \right) \le  \exp \left(
   2 \gamma_n \frac{R}{1-R} \right)
\end{equation*}
for $z \in \D_R:=\{z \in \C:|z|<R\}$. We use $f(0)=0$ once more to apply the
Schwarz' lemma for
the disk $\D_R$ and obtain
\begin{equation}\label{eq:f}
 |f(z)| \le \frac{1}{R} \exp \left( 2 \gamma_n \frac{R}{1-R} \right) \,
|z| \, , \qquad z\in \D_R \, . 
\end{equation}
We next consider the function
$$R \mapsto  \frac{1}{R} \exp \left( 2 \gamma_n \frac{R}{1-R} \right)\, , \qquad R
\in (0,1)\,. $$  
It is easy to check that this function attains its minimal value on $
(0,1)$ and there is exactly one such point. Call it $R_n$. Then a straightforward
computation gives  formula (\ref{eq:Rn}) for $R_n$. Inserting (\ref{eq:Rn})
into (\ref{eq:f}) leads to (\ref{eq:schwarz}).
\hfill{$\blacksquare$}

\bigskip

{\bf Proof of Remark \ref{rem:2}.}
First, $(R_n)$ is a strictly increasing sequence which tends to the limit
$1$ because the function $x \mapsto
1+x-\sqrt{x^2+2x}$ is strictly decreasing on $[0,\infty)$ and the sequence
$(\gamma_n)$ is strictly decreasing with limit $0$. Once again, since   $(\gamma_n)$ is
strictly decreasing, the expression
\begin{equation*}
\frac{1}{R} \exp \left( 2 \gamma_n\, \frac{R}{1-R} \right)\, , \qquad R\in (0,1)\,,
\end{equation*}
is strictly decreasing with respect to $n$, so its minimal value
\begin{equation*}
\frac{1}{R_n} \exp \left( 2 \gamma_n \frac{R_n}{1-R_n} \right)=\frac{1}{R_n} \exp \left( \sqrt{\gamma_n^2+\gamma_n}-\gamma_n \right)
\end{equation*}
 is also strictly decreasing with respect to $n$ (see the proof of Theorem \ref{thm:2a}).

Consequently, 
\begin{equation*}
\left( \frac{1}{R_n} \exp \left( \sqrt{\gamma_n^2+\gamma_n}-\gamma_n \right) \right)
\end{equation*}
forms  a strictly decreasing sequence with limit $1$.

\medskip

Now, let $f \in \mathcal{A}_n$. Then inequality (\ref{eq:schwarz}) implies
\begin{equation*}
|f'(0)| \le \frac{1}{R_n} \exp \left(  \sqrt{\gamma_n^2+2 \gamma_n}-\gamma_n\
\right)\,.
\end{equation*} 
A straightforward but lengthy calculation, using Theorem \ref{thm:gamma_n} (a), gives 
\begin{equation*}  
|f'(0)| \le \frac{1}{R_n} \exp \left(  \sqrt{\gamma_n^2+2 \gamma_n}-\gamma_n\
\right)=1+8 \frac{\Gamma(5/4)}{\Gamma(3/4)} \frac{1}{\sqrt{n}}+ O(1/n) \qquad
(n \to \infty)\,.
\end{equation*}

If $f_n \in \mathcal{A}_n$ is a universal covering map of $\C\backslash S_n$,
then $f_n^*\lambda_{\C \backslash S_n}(z) \, |dz|=\lambda_{\D}(z)\, |dz|$. Therefore,  
\begin{equation*}
|f_n'(0)|=\frac{2}{\lambda_{\C\backslash S_n}(0)}=\lim \limits_{z \to 0}
\frac{2}{n |z|^{n-1} \lambda_{1-1/n,1,1}(z^n)}
\end{equation*}
by  Theorem \ref{thm:3a}. Combining this with  Theorem \ref{thm:formula} and
(\ref{eq:K3}) as well as (\ref{eq:K2/K3})
we get

\begin{equation*}
\begin{split}
 |f_n'(0)|= \frac{1}{n} &\left(
   -\frac{\sin\left(\frac{\pi}{n}\right)}{\pi}\,\, 
\Gamma \left(\frac{n-1}{2   n}\right)^2 
\Gamma \left(\frac{n+1}{2 n}\right)^2 
\, \, _2F_1\left( \frac{n-1}{2 n},\frac{n-1}{2 n},1;1\right)^2 \right.
\\[2mm]
&\left. \quad + \,  2 \, \frac{\sin\left(\frac{\pi}{n}\right)}{\pi}\, \,
\Gamma \left(\frac{n-1}{2
   n}\right)^2 \Gamma\left(\frac{1}{n} \right)\,\, _2F_1\left(\frac{n-1}{2 n},\frac{n-1}{2 n},1;1\right)
  \right)\,.
\end{split}
\end{equation*}

Because of Gauss's theorem 
\begin{equation*}
_2F_1(a,b,c;1)=\frac{\Gamma(c)
  \Gamma(c-a-b)}{\Gamma(c-a) \Gamma(c-b)} \quad \text{ for } \quad \Re (c)>\Re
(a+b),\, \, 
c\not=0,-1,-2,\ldots\,,
\end{equation*}
  (see \cite[p.~213]{AS84}), we get
\begin{equation*}
 |f_n'(0)|=\frac{\sin \left(\frac{\pi }{n}\right) \Gamma \left(\frac{1}{n}\right)^2 \Gamma \left(\frac{n-1}{2 n}\right)^2}{\pi  n \Gamma
   \left(\frac{n+1}{2 n}\right)^2}\, .
\end{equation*}

This simplifies, using the reflection formula (\ref{eq:reflection}), to
\begin{equation*}
 |f_n'(0)|=\frac{\Gamma \left(\frac{1}{2} \left(1-\frac{1}{n}\right)\right)^2 \Gamma
   \left(\frac{1}{n}\right)}{n \Gamma \left(1-\frac{1}{n}\right) \Gamma \left(\frac{1}{2}
   \left(1+\frac{1}{n}\right)\right)^2} \, .
\end{equation*}
Using the well--known Taylor series expansion 
\begin{equation*}
1/\Gamma(z)= z + \gamma z^2 +\cdots
\end{equation*}
 at $z=0$, where $\gamma$ is the Euler--Mascheroni constant, (cf.~\cite[p.~76/77]{AS84}), one readily
verifies that
\begin{equation*}
 \frac{\Gamma \left(\frac{1}{2} \left(1-\frac{1}{n}\right)\right)^2 \Gamma
   \left(\frac{1}{n}\right)}{n \Gamma \left(1-\frac{1}{n}\right) \Gamma \left(\frac{1}{2}
   \left(1+\frac{1}{n}\right)\right)^2} =1+\frac{4 \log 2}{n}+O(1/n^2) \qquad
(n \to \infty) \, . 
\end{equation*}
The proof of Remark \ref{rem:2} is now complete.
\hfill{$\blacksquare$}

\begin{figure}[h]
\centerline{\includegraphics{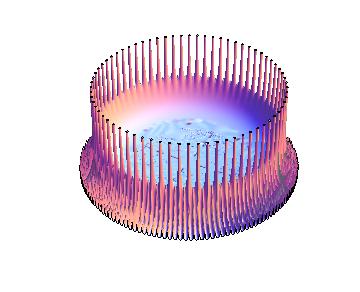}}

\vspace*{-1.5cm}
\caption{The hyperbolic density $\lambda_{\C\backslash S_n}$ for $n=70$}
\end{figure}

\vspace*{1cm}

Daniela Kraus\\ 
Department of Mathematics\\
University of W\"urzburg\\
Emil Fischer Stra{\ss}e 40\\
97074 W\"urzburg\\
dakraus@mathematik.uni-wuerzburg.de\\
Phone: +49 931 318 5028

\vspace{0.5cm}

Oliver Roth\\
Department of Mathematics\\
University of W\"urzburg\\
Emil Fischer Stra{\ss}e 40\\
97074 W\"urzburg\\
roth@mathematik.uni-wuerzburg.de\\
Phone: +49 931 318 4974


\begin{thebibliography}{99}
\bibitem{AS84} M.~Abramowitz and I.A.~Stegun, {\it Pocketbook of Mathematical
    Functions}, Harri Deutsch, Frankfurt, 1984.

\bibitem{Ahl38} L.~Ahlfors,  An extension of Schwarz's
lemma, {\it Trans.~Amer.~Math.~Soc}.~(1938), {\bf 43}, 359--364.


\bibitem{Dra69} D.~Drasin, Normal families and the Nevanlinna theory,
  \textit{Acta Math}.~(1969), \textbf{122}, 231--263.

\bibitem{GT97} D.~Gilbarg and N.~S.~Trudinger, {\it Elliptic Partial
    Differential Equations of Second Order}, Springer, Berlin--New York, 1997.

\bibitem{Hay} W.~K.~Hayman, Some remarks on Schottky's theorem, {\it
    Proc.~Cambrdige Philos.~Soc}. (1947), {\bf 43}, 442--454.


\bibitem{Hei62} M.~Heins, 
On a class of conformal metrics, {\it Nagoya Math.~J}.~(1962), {\bf 21},
    1--60.


\bibitem{Hem79}
J.~A.~Hempel, 
The Poincar\'e metric on the twice punctured plane and the theorems of 
Landau and Schottky,
{\it J.~Lond.~Math.~Soc., II.~Ser}.~(1979), {\bf  20}, 435--445.

\bibitem{Hem80}
J.~A.~Hempel, Precise bounds on the theorems of Schottky and Picard,
{\it J.~London Math. Soc}.~(1980), {\bf 21}, 279--286.

\bibitem{Jen} J.~Jenkins, On explicit bounds in Landau's theorem II, {\it
    Can.~J.~Math}.~(1981), {\bf 33}, 559-562.

\bibitem{KR2007} D.~Kraus and O.~Roth, The behaviour of solutions of the  Gaussian curvature
equation near an isolated boundary point, 
{\it Math.~Proc.~Cambr.~Phil.~Soc}.~(2008), {\bf 145}, 643--667.

\bibitem{KRS} D.~Kraus, O.~Roth and T.~Sugawa,
 Metrics with conical singularities on the  sphere
 and sharp extensions of the theorems of Landau and Schottky,
 \textit{Math.~Z}.~(2011), \textbf{267} No.~3--4,
 851--868.

\bibitem{Lan04} E.~Landau, \"Uber eine Verallgemeinerung des Picardschen
  Satzes, \textit{S.B.~Preuss.~Akad. Wiss.}~(1904), 1118-1133.

\bibitem{LQ2007} Z.~Li and Y.~Qi, A remark on Schottky's theorem,
{\it Bull.~London Math.~Soc}.~(2007), {\bf 39}, 242--246.

\bibitem{Min87} D.~Minda, The strong form of Ahlfors' lemma, {\it Rocky
    Mountain J.~Math}.~(1987), {\bf 17} no.~3, 457--461.

\bibitem{Min97} D.~Minda, The density of the hyperbolic metric near an
  isolated boundary point, {\it Complex Variables}~(1997), {\bf 32}, 331--340.

\bibitem{Mon27} P.~Montel, \textit{Le\c{c}ons sur les familles normales des
    fonctions analytiques et leurs applications}, Gauthier-Villars, Paris, 1927.


\bibitem{Nit57} J.~Nitsche, \"Uber die isolierten Singularit\"aten der
  L\"osungen von $\Delta u =e^u$, {\it Math.~Z.}~(1957), {\bf 68}, 316--324.


\bibitem{Ran95} T.~Ransford, {\it Potential Theory in the Complex Plane},
  Cambridge University Press, Cambridge, 1995.

\bibitem{Roy86} H.~L.~Royden, The Ahlfors--Schwarz lemma: the
  case of equality, {\it J.~Analyse Math}.~(1986), {\bf 46}, 261--270.

\bibitem{Sch04} F.~Schottky,  \"Uber den Picardschen Satz und die Borelschen
  Ungleichungen, \textit{S.B.~Preuss. Akad.~Wiss.},  1244--1263, 1904. 

\bibitem{Zal98} L.~Zalcman, Normal families: New perspectives,
\textit{Bull.~Amer.~Math.~Soc}.~(1998), \textbf{35} No.~3, 215-230.
\end{thebibliography}
\end{document}